\documentclass[10pt]{article}
\usepackage[affil-it]{authblk}

\usepackage{amsmath, amsthm, amssymb}
\usepackage{amssymb}
\usepackage{graphicx, color}
\usepackage{subfigure}
\usepackage{url}

\newcommand\bi{\begin{itemize}}
\newcommand\ei{\end{itemize}}

\def\addsymbol #1: #2#3{$#1$ \> \parbox{5in}{#2 \dotfill \pageref{#3}}\\}

\newtheorem{defin}{Definition}[section]
          
          \newtheorem{teo}[defin]{Theorem}
          \newtheorem{con}[defin]{Conjecture}
          \newtheorem{cond}[defin]{Condition}
          \newtheorem{prop}[defin]{Proposition}
          \newtheorem{lem}[defin]{Lemma}
          \newtheorem{rmk}[defin]{Remark}
          \newtheorem{cor}[defin]{Corollary}

          \newcommand{\bfact}{\begin{fact}}
          \newcommand{\efact}{\end{fact}}

          \newcommand{\beq}{\begin{equation}}
          \newcommand{\eeq}{\end{equation}}
          \newcommand{\beqn}{\begin{eqnarray}}
          \newcommand{\beqnn}{\begin{eqnarray*}}
          \newcommand{\eeqn}{\end{eqnarray}}
          \newcommand{\eeqnn}{\end{eqnarray*}}
          \newcommand{\bprop}{\begin{prop}}
          \newcommand{\eprop}{\end{prop}}
          
          \newcommand{\bc}{\be\begin{array}{r@{\,}c@{\,}l}}
\newcommand{\ec}{\end{array}\ee}

          \newcommand{\bcor}{\begin{cor}}
          
          \newcommand{\ecor}{\end{cor}}
          \newcommand{\bcon}{\begin{con}}
          \newcommand{\econ}{\end{con}}
          \newcommand{\bcond}{\begin{cond}}
          
          \newcommand{\econd}{\end{cond}}
          \newcommand{\bteo}{\begin{teo}}
          \newcommand{\eteo}{\end{teo}}
          \newcommand{\brm}{\begin{rmk}}
          \newcommand{\erm}{\end{rmk}}
          \newcommand{\blem}{\begin{lem}}
          \newcommand{\elem}{\end{lem}}

          \newcommand{\ben}{\begin{enumerate}}
          \newcommand{\een}{\end{enumerate}}
          \newcommand{\bei}{\begin{itemize}}
          \newcommand{\eei}{\end{itemize}}
	 
          \newcommand{\bdf}{\begin{defin}}
          \newcommand{\edf}{\end{defin}}

          \renewcommand{\>}{&>&}

          \newcommand{\Z}{{\mathbb Z}}

          \newcommand{\R}{{\mathbb R}}
          
          \newcommand{\E}{{\mathbb E}}
          
          \renewcommand{\P}{{\mathbb P}}

          \newcommand{\N}{{\mathbb N}}

          \newcommand{\W}{{\cal W}}
          
          \newcommand{\cH}{\cal H}

          \newcommand{\om}{\omega}

 \newcommand\epl{{\eps,+}}
 \newcommand\emi{{\eps,-}}

\newcommand{\btt}{\begin{theorem}}
\newcommand{\ett}{\end{theorem}}

\newcommand{\Wl}{\mathcal{W}_l}
\newcommand{\Wr}{\mathcal{W}_r}
\newcommand{\be}{\begin{equation}}
\newcommand{\ee}{\end{equation}}

\newcommand\eps{\epsilon}

\newcommand\epm{{\eps,\pm}}

          \newcommand\sqr{\vcenter{
          \hrule height.1mm
          \hbox{\vrule width.1mm height2.2mm\kern2.18mm\vrule width.1mm}
          \hrule height.1mm}}        

\newcommand{\astoo}[2]{\underset{{#1}\to{#2}}{\longrightarrow}}

\newcommand{\Astoo}[1]{\underset{{#1}\to 0}{\Longrightarrow}}

\begin{document}

\title{Perturbations of supercritical oriented percolation  \\
and sticky Brownian webs}
\author{Emmanuel Schertzer$^{\,1,2}$ \and Rongfeng Sun$^{\,3}$}
\maketitle

\footnotetext[1]{Laboratoire de Probabilit\'es, Statistiques et Mod\'elisation (LPSM), Sorbonne Universit\'e,
CNRS UMR 8001, Paris, France. Email: emmanuel.schertzer@upmc.fr}

\footnotetext[2]{
Center for Interdisciplinary Research in Biology (CIRB), Coll\`ege de France, CNRS UMR 7241,
PSL Research University, Paris, France
}

\footnotetext[3]{Dept.\ of Mathematics, National University of Singapore, 10 Lower Kent Ridge Road, 119076 Singapore. Email: matsr@nus.edu.sg}

\centerline{\em To our advisor Chuck Newman, on his 70th birthday}

\begin{abstract}
Previously, Sarkar and Sun \cite{SS13} have shown that for supercritical oriented percolation in dimension $1+1$, the set of rightmost infinite open paths converges to the Brownian web after proper centering and scaling. In this note, we show that a pair of sticky Brownian webs arise naturally if one considers the set of right-most infinite open paths for two coupled percolation configurations with distinct (but close) percolation parameters.
This leads to a natural conjecture on the convergence of the dynamical supercritical oriented percolation model to the so-called dynamical Brownian web.
\end{abstract}

\noindent
{\it AMS 2010 subject classification:} 60K35, 82B43.\\
{\it Keywords.} Brownian net, Brownian web, oriented percolation, sticky Brownian motion.
\vspace{12pt}

\section{Main result}\label{intro}

\noindent

{\bf Oriented Percolation and the Brownian Web.} Following the notation of \cite{SS13}, we consider supercritical oriented percolation on $\Z^2_{\rm even}:=\{(x,t)\in \Z^2: x+t \mbox{ is even}\}$, with directed edges leading from $(x,t)$ to $(x\pm 1, t+1)$ for each $(x,t)\in\Z^2_{\rm even}$. Independently, each directed edge is open with probability $p$ and closed with probability $1-p$ for some $p>p_c$, the critical percolation parameter. A site $z\in\Z^2_{\rm even}$ is said to be a percolation point if there is an infinite open path of directed edges starting from $z$, and let ${\cal K}\subset \Z^2_{\rm even}$ denote the set of percolation points. From each $z=(x,t)\in {\cal K}$, let $\gamma_z :=(\gamma_z(n))_{n\geq t}$ denote the rightmost infinite open path starting from $z$. When $z=(x,t)\not \in {\cal K}$, define $\gamma_z$ to be the rightmost infinite open path starting from $(-\infty, x]\times\{t\}$.

\begin{figure}
  \includegraphics[width=\linewidth]{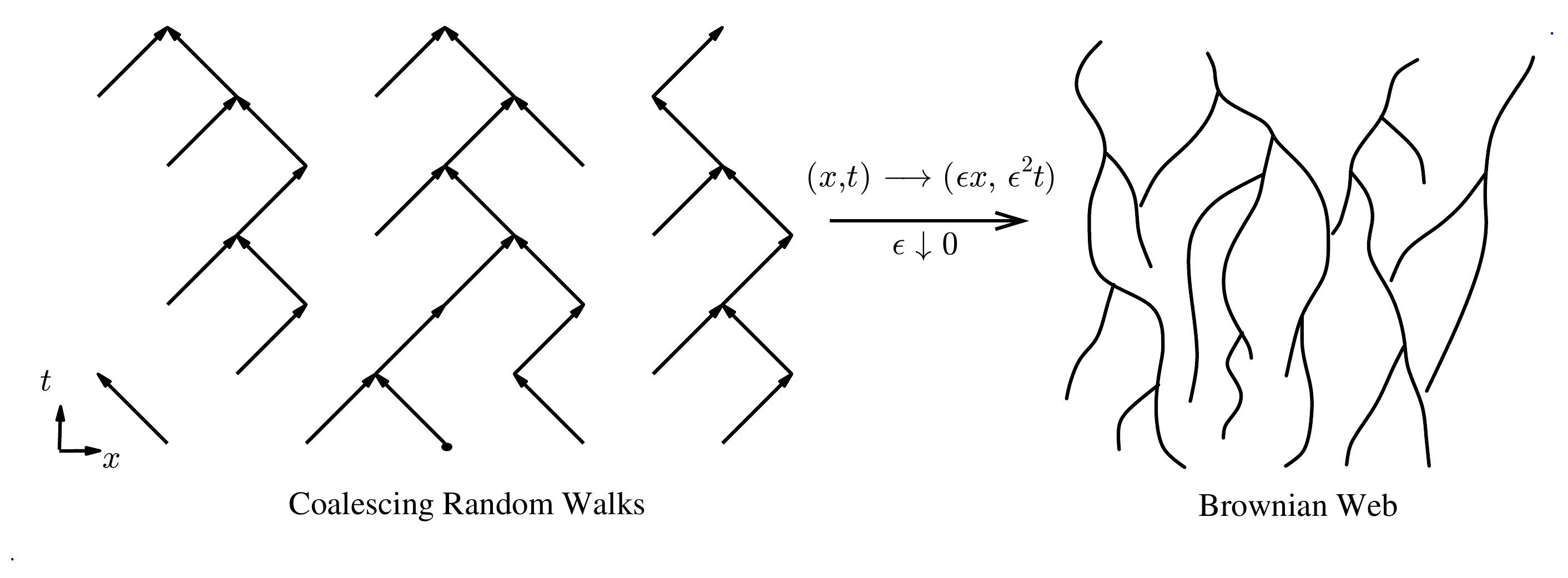}
  \caption{The Brownian web as scaling limit of coalescing random walks. }
  \label{fig:web}
\end{figure}

In \cite{SS13}, Sun and Sarkar showed that after suitable centering and rescaling, the set of paths $ \Gamma=\{\gamma_z \ : \  z\in{\cal K} \}$ converges to a family of coalescing Brownian motions, known as the Brownian web \cite{FINR04, TW98, SSS17} and which can  also be obtained as the ``universal" scaling limit  of coalescing random walks \cite{NRS05}. See Fig \ref{fig:web}. More precisely, for $a, b, \eps>0$, let
\begin{equation}\label{Seps}
S_{a, b, \eps}(x,t) := \Big(\frac{\eps}{b}(x-at), \eps^2 t\Big), \qquad (x,t)\in \R^2,
\end{equation}
be a space-time scaling map, and let $S_{a, b, \eps}(\gamma_z)$ and $S_{a,b,\eps}(\Gamma)$ be defined by identifying each path with its graph in $\R^2$. Then for $p>p_c$, there exist $\alpha(p), \sigma(p)>0$ such that:
\begin{equation}\label{cv:web-ss}
S_{\alpha(p),\sigma(p),\eps}(\Gamma) \Astoo{\eps} {\cal W},
\end{equation}
where $\Rightarrow$ denotes weak convergence, $\cal W$ is the {\em Brownian web}, and both $S_{\alpha(p),\sigma(p),\eps}(\Gamma)$ and $\cal W$ are regarded as random variables taking values in $\cH$, a suitable space of compact sets of paths. See \cite{SS13} or the recent survey \cite{SSS17} for further details details on the Brownian web, including the space $\cal H$ in which the Brownian web takes its value.

{\bf  Perturbed percolation clusters and sticky Brownian webs.} Let $p>p_c$ and let $\Omega$ be the percolation configuration with percolation parameter $p$.  We now describe two natural perturbations of
$\Omega$  which consists in slightly increasing (resp., decreasing)
the fraction of open edges.

In order to do so, we equip each directed edge $e$ of $\Z_{\rm even}^2$ with an independent uniform random variable $\om_e$ on $[0,1]$. For any $p
\in[0,1]$,  we declare an edge $e$ to be $p$-open iff $p<\om_e$. This construction provides a natural coupling between percolation configurations with different percolation parameters $p$. For $\eps<(p-p_c)\wedge (1-p)$, we consider the percolation configurations $\Omega^{\eps,-}$ and $\Omega^{\eps,+}$
corresponding to the parameter $p-\eps$ and $p+\eps$ respectively, and which are coupled in such a way so that $\Omega^{\eps,+}$ has a surplus of open edges.
Let ${\cal K}^{\eps,\pm}\subset \Z^2_{\rm even}$ denote the set of percolation points in the configuration $\Omega^{\eps,\pm}$. For every $z\in{\cal K}^{\eps,\pm}$, let $\gamma^{\eps,\pm}_z$ denote the rightmost infinite open path starting from $z$, and let $\Gamma^{\eps, \pm}$ denote the set of all such paths. We note that for $\eps$ small enough, ${\cal K}^{\eps,\pm}$ is non-empty since $p>p_c$.

\medskip

The first aim of this note is to describe the correlation between the sets of right-most paths $\Gamma^{\eps,+}$ and $\Gamma^{\eps,-}$, thus extending the link between oriented percolation and the Brownian web established in \cite{SS13}. In order to explain our results, we recall the definition of  left-right sticky  pair of Brownian webs. This definition is slightly different from the one in Theorem 1.5 \cite{SS08}
(where this object first appeared) and is adapted from the definition of $\theta$-coupled webs in Theorem 4 in \cite{HW09} or Theorem 76 in \cite{H07}.
We will comment more on the connection between the definitions in \cite{SS08} and \cite{HW09} in the Appendix. We first recall the following from \cite[Proposition~2.1]{SS08}.

\begin{defin}[Left-right sticky pair of Brownian paths]\label{def:lrsde}
Let $(l,r)$ be two paths with repective starting time $\sigma_l, \sigma_r$.
Assume  first that  $\sigma = \sigma_l=\sigma_r$. Then $(l,r)$
is distributed as a left-right sticky pair  of Brownian paths with drift $b$ iff
the pair $(l, r)$ is distributed as the unique weak solution of the {\em left-right SDE}
\begin{equation}\label{lrsde}
\begin{aligned}
\forall\, t\geq \sigma, \ \ \  d \bar l_t &= 1_{\bar l_t\neq \bar r_t} dB^l_t + 1_{\bar l_t= \bar r_t} dB^s_t - bdt, \\
d \bar r_t &= 1_{\bar l_t\neq \bar r_t} dB^r_t + 1_{\bar l_t= \bar r_t} dB^s_t + bdt,
\end{aligned}
\end{equation}
(where $B^l, B^r, B^s$ are independent standard Brownian motions) with initial condition $(l(\sigma), r(\sigma))$ at time $\sigma$, and
subject to the constraint
\begin{enumerate}
\item[$(${\bf C}$)$]\label{cond:c} $\bar l(t)\leq \bar r(t)$ for all $t\geq \tau:=\inf\{s: \bar l(s)\leq \bar r(s)\}.$
\end{enumerate}
If $\sigma_l \neq \sigma_r$, then $(l,r)$ will be called a left-right sticky pair if conditioned on $(l_t, r_t)_{t\leq \sigma_l\vee \sigma_r}$,
$(l_t, r_t)_{t\geq \sigma_l\vee\sigma_r}$ solves the left-right SDE \eqref{lrsde}.
\end{defin}

\brm The term {\em sticky} is motivated by the fact that the process $w=\frac{l-r}{\sqrt{2}}$ is a drifted Brownian motion with sticky reflection at $0$ (see e.g. \cite{W02}) when $l$ and $r$ have the same starting point.
\erm

The following theorem is adapted from \cite[Theorem 4]{HW09} on the characterization of a pair of sticky Brownian webs with zero drift. We will sketch the proof in the Appendix, where we show its equivalence with the characterization given in \cite[Theorem 1.5]{SS08}.

\begin{teo}[Left-right sticky Brownian webs]\label{sticky-webs}
In distribution, there exists a unique pair $(\Wl, \Wr)$ valued in  $\cH\otimes\cH$ such that
\begin{enumerate}
\item $\Wl$ (resp., $\Wr$) is a Brownian web with drift $-b$ (resp., $b$). In particular, for every deterministic  $z$, there is a.s.\ a unique path $l_{z}$ (resp., $r_z$) in $\Wl$ (resp., $\Wr$)  starting from $z$.
\item For every deterministic pair $z_1,z_2\in\R^2$, the pair $(l_{z_1},r_{z_2})\in\Wl\otimes\Wr$ is distributed as a left-right sticky pair of Brownian paths with drift $b$.
\item[3](Co-adaptedness) Let $({\cal G}_t; t\in\R)$ be the natural filtration induced by the pair of Brownian webs $(\Wl,\Wr)$ (see Definition \ref{def:common-filtration} below for more details). For any deterministic $z_1,z_2\in \R^2$, the pairs of processes
$(r_{z_1},r_{z_2}), (l_{z_1}, l_{z_2})$ and $(r_{z_1},l_{z_2})$
are all Markov with respect to $({\mathcal G}_t; t\geq0)$, i.e., for any such pair $(l, r)$, $P( (l(s), r(s))_{s\geq t}\in \cdot |{\cal G}_t) = P( (l(s), r(s))_{s\geq t}\in \cdot  |l(t), r(t))$.
\end{enumerate}
We call $(\Wl, \Wr)$ a left-right  pair of sticky Brownian webs with drift $b>0$.
\end{teo}

%

Our main result is the following, which generalizes the main result of \cite{SS13} from a single Brownian web to a pair of sticky Brownian webs. The constants $\sigma(p)$ and $\alpha(p)$ are as in
the convergence statement (\ref{cv:web-ss}) and as in \cite[(2.4)]{SS13}, and were first introduced in \cite{K89}. See also (\ref{def-sigma-alpha}) below
for more details.

\bteo\label{thm:main}

The function $p\to \alpha(p)$ is differentiable and $\sigma(p)>0$ on $(p_c,1]$ . Furthermore,
$$
S_{\sigma(p),\alpha(p),\eps}(
\Gamma^{\eps,-}, \Gamma^{\eps,+} ) \Astoo{\eps} (\W^l, \W^r),
$$
where $ (\W^l, \W^r)$ is a left-right pair of sticky Brownian webs with drift $b:=\frac{\alpha'(p)}{\sigma(p)}$, regarded as a random variable taking values in $\cH\otimes\cH$ (see \cite{SSS17} for further details).
\end{teo}

\begin{rmk}
We note that a straightforward extension of \cite{SS13} will show the convergence of
the marginals, namely
$S_{\sigma(p),\alpha(p),\eps}(\Gamma^{\eps,-} ) \Astoo{\eps} \Wl$ and
$S_{\sigma(p),\alpha(p),\eps}(\Gamma^{\eps,+} ) \Astoo{\eps} \Wr$, so most of our work will be dedicated to showing that the two webs have a sticky interaction (conditions 2-3 of Theorem \ref{sticky-webs}).
\end{rmk}




\bigskip

\noindent
{\bf Discrete web perturbations and the dynamical Brownian web.} We close the introduction with a natural conjecture arising from the previous result. In order to do so, we make a detour, and we consider an alternative (and simpler) model where sticky Brownian webs also emerge naturally. Independently at each edge, declare one of the  two randomly chosen out-going edge (or arrow) $(x\pm1,t+1)$ to be open, whereas the other edge remains closed. At any point  $z\in\Z_{even}^2$, there is exactly a single (infinite) path $w_z$ starting from $z$ and we denote by $W$ the infinite collection of paths $\{w_z\}_{z\in \Z_{even}^2}$. This set is often referred to as the {\it discrete web} that can be described as  an infinite set of coalescing random walks. In \cite{FINR04}, it was shown that
\[S_{0,1,\eps}(W) \Astoo{\eps}  \W,\]
where $\W$ is the Brownian web.
Let us now perturb the previous arrow configuration by adding an additional arrow independently at each site $z\in\Z_{even}^2$ with probability $\eps\in[0,1]$. For $\eps>0$, at any point $z$, there are infinitely many outgoing paths due to the branching of arrows, and we define $\tilde \Gamma^{+,
 \eps}$ (resp.,  $\tilde \Gamma^{-,
 \eps}$) to be the the infinite collection  of right-most (resp., left-most paths) starting from $\Z_{even}^2$. In \cite{SS08}, it was shown that
 $$
 S_{0,1,\eps}(\tilde \Gamma^{-,\eps},\tilde\Gamma^{+,\eps}) \Astoo{\eps}  \ (\Wl,\Wr)
 $$
 where $(\Wl,\Wr)$ is a pair left-right sticky Brownian webs with drift $1$. In light of Theorem \ref{thm:main}, we observe that the left-right sticky pair of webs arises as a perturbation of the underlying path configuration in two distinct models: the discrete  web and the supercritical percolation model.

\medskip

In the discrete web described above, the left-right sticky pair of Brownian webs was obtained by adding some extra branching point in the underlying arrow configuration. We now describe an alternative  perturbation of the discrete web $W$  where sticky webs also arise naturally.

Start with the discrete web $W$ and equip each vertex  $z\in\Z_{\rm even}^2$ with
an independent Poisson clock with rate $\eps$. Every time the clock rings at a given vertex $z$,  switch the direction of the arrow starting from $z$ (i.e., the edge $((x,t),(x\pm1))$ becomes
$((x,t),(x\mp1))$). This defines a stationary Markov process
$(W_s^
\eps; s\geq0)$ -- called the discrete dynamical web  \cite{FNRS09} --  whose one-dimensional marginal is given by the law of the discrete web. In \cite{NRS10}, it was shown that
\[S_{0,1,\eps}(W_s; s\geq0) \ \Astoo{\eps}  (\W_s; s\geq0)\]
where $(\W_s; s\geq0)$ is a continuum objet called the dynamical Brownian web and the convergence is meant in the sense of finite dimensional distribution. See also \cite{HW09} for a definition of the continuum object. Not surprisingly, this process is a stationary Markov process with marginal distribution being  the law of the Brownian web, and the two-dimensional  marginal $(\W_{s_1}, \W_{s_2})$ can be described in terms  of a pair of sticky Brownian webs.  In this context, the definition of a sticky pair of webs is analogous to Definition
\ref{sticky-webs} with the notable difference that the Brownian motions have no drift and there is no natural ordering of the webs $(\W_{s_1}, \W_{s_2})$ (see condition (C) in Deinition \ref{sticky-webs}). See \cite{HW09,NRS10} for more details.

\bigskip

\noindent

{\bf Dynamical percolation.} Let us now re-consider the supercritical oriented percolation model with parameter $p>p_c$. Analogously to the dynamical web, we can equip every edge with an independent continuous-time Markov chain $(v_z(s); s\geq0)$ with initial law $\P(v_e(0)=1)=p$ and $\P(v_e(0)=0)=1-p$, and the chain switches from state $0$ to $1$ at rate $\eps p$ and switches from $1$ to $0$ at rate $\eps(1-p)$. Note that the initial law of the Markov chain is also stationary under the dynamics.
For every ``dynamical time'' $s\geq 0$, we declare an edge to be open if and only if $v_e(s)=1$. This defines a stationary dynamical oriented percolation model, whose marginal distribution at each dynamical time is the (static) oriented percolation model on $\Z^2_{\rm even}$ with parameter $p$.

\begin{con}
Let $\Gamma_s^\eps$ be the set of right-most infinite open paths at dynamical time $s$.
There exists a constant $c(p)$ such that
$$S_{\alpha(p),\sigma(p),\eps}(\Gamma_{s}^\eps; s\geq0) \Astoo{\eps}   (\W_{c(p)s}; s\geq0)$$
where $(\W_s; s\geq0)$ is the dynamical Brownian web.
\end{con}

We note that scaling limits of dynamical  dynamical percolation models have been obtained in the case of standard percolation in $\Z^2$ \cite{CFN06, GPS13, GPS13bis}, but only at criticality (i.e. when $p=p_c$). In the context of oriented percolation, the previous conjecture would show that some interesting large-scale dynamics takes place even in the super-critical regime.

\bigskip

{\bf Technical remarks.}
As pointed out before Theorem \ref{sticky-webs}, the characterization of a sticky pair of left-right Brownian webs presented in this paper is not new, and is adapted from the one of $\theta$-coupled Brownian webs given in \cite[Theorem 4]{HW09}. However, to the best of our knowledge, this characterization has never been used to prove scaling limit results, and we believe that it could be of interest in other settings. The latter characterization is particularly efficient since convergence boils down to proving that (1) $\Wl$ and $\Wr$ are two Brownian webs, and (2) each pair of paths with deterministic starting points has the required distribution w.r.t. to the {\em larger filtration $({\cal G}_t; t\geq0)$}. We will discuss more about this in the appendix.

Finally, we would also like to draw attention to Proposition \ref{prop:lr}, which  provides a minimal characterization of a left-right sticky pair of Brownian paths that could be useful for proving convergence in other contexts.

\section{Preliminaries}\label{sect:preliminaries}
\begin{figure}
\centering
  \includegraphics[scale=0.4]{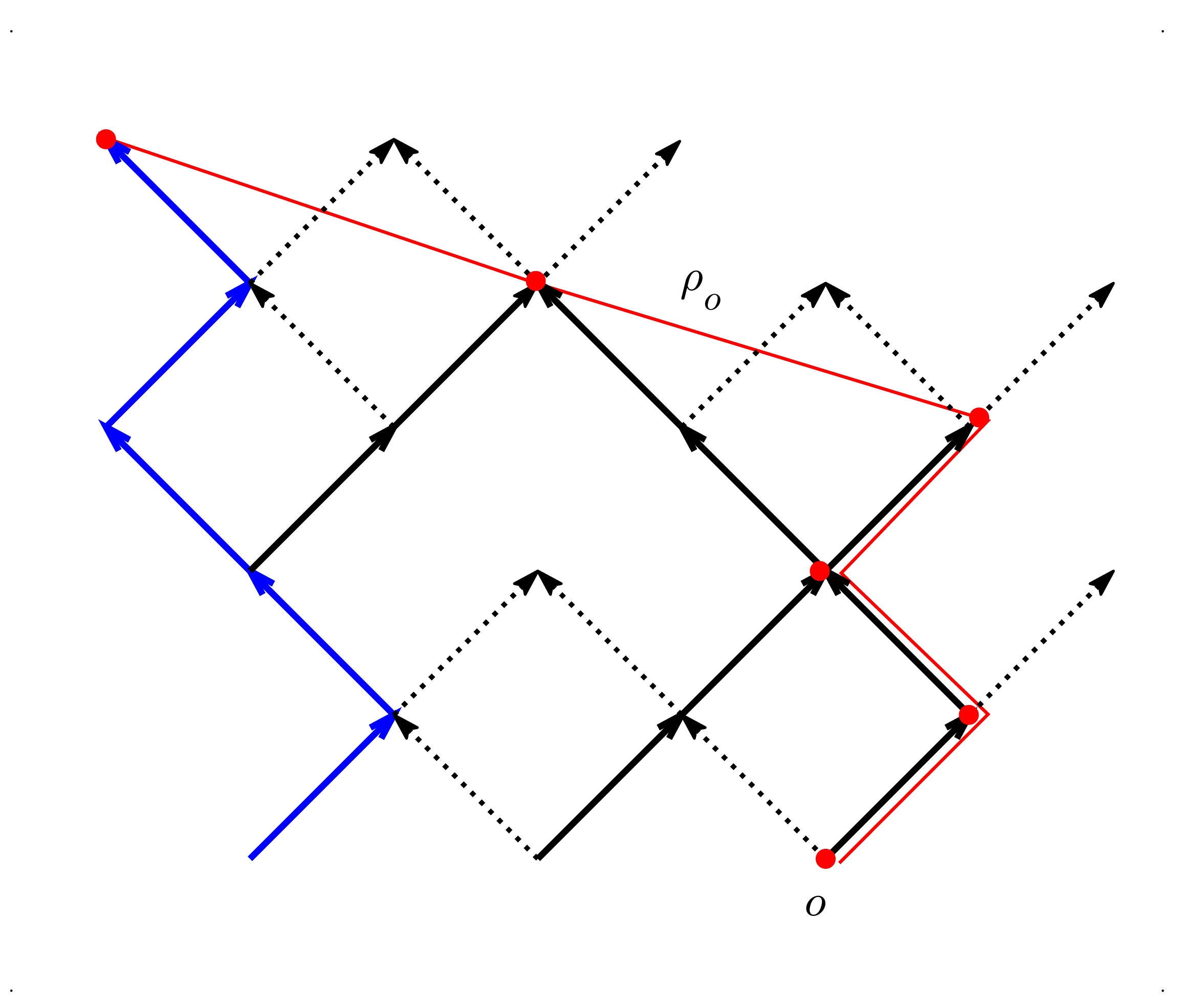}
  \caption{The path $\rho_o$ and the exploration cluster starting from $o$. $\rho_o(n)$ is the rightmost position at time $n$ that can be reached by some open path starting from $(-\infty, 0]\times\{0\}$. The exploration cluster associated to $o$ is obtained by exploring the set of vertices  connected to $(-\infty, 0]\times\{0\}$ from right to left. See \cite{SS13} for a precise description of the exploration algorithm to generate the cluster $C_{o}$.}
  \label{fig:explore}
\end{figure}

The main tool used in \cite{SS13} is the notion of a {\em percolation exploration cluster}, which we briefly recall here. For every $z=(x,t)\in\Z_{\rm even}^2$, the percolation exploration cluster $C_z:=(C_z(n))_{n\geq t}$ consists of a set of sequentially explored edges such that for each time $n>t$, a minimal set of edges $C_z(n)$ before time $n$ are explored in order to find the rightmost open path connecting $(\infty, x]\times\{t\}$ to $\Z\times\{n\}$. See Figure \ref{fig:explore}.
The percolation exploration cluster $C_z$ provides a good approximation for the rightmost infinite open path $\gamma_z$. Let $\rho_z$ be the analog of $r_z$ in \cite{SS13}, i.e., $\rho_z(n)$ is the rightmost position at time $n$ that can be reached by some open path starting from $(-\infty, x]\times\{t\}$. It was shown in \cite{SS13} that $C_z$ is bounded between the paths $\gamma_z$ and $\rho_z$, and furthermore, $\gamma_z$ and $\rho_z$ converge to the same Brownian motion after proper centering and scaling. The advantage of approximating $\gamma_z$ by $\rho_z$ is that the latter depends only on the edge configurations explored up to the present and not on the future.

As in \cite{SS13},  we denote by $(\rho(T_i),T_i)_{i\in\N}$ the successive break points along $\rho:=\rho_o$. More precisely, the $T_i$'s correspond to the successive times at which $\gamma:=\gamma_o$ and $\rho$ coincide. Let
\begin{align*}
& \tau_1 := T_1, \qquad  &X_1 &:= \rho(T_1), \\
& \tau_i := T_i - T_{i-1}, &X_i &:= \rho(T_i) - \rho(T_{i-1})  \qquad \mbox{for $i\geq 2$.}
\end{align*}
It is known from \cite{K89} that $((X_i,\tau_i); \ i\geq 2)$ forms a sequence of i.i.d. random variables
with $\tau_2$ and $X_2$ having all finite moments. We will prove the following result.

\bprop\label{lem:f-g-h-diff}
For every $i,j\in\N$, the function
\begin{equation}\label{fijp}
f_{i,j}(p):= \E_p[X_2^i \tau_2^j]
\end{equation}
is differentiable on $(p_c,1)$.
\eprop

As a corollary, we have
\bcor\label{cor:diif-alpha-sigma}
$p\to \alpha(p)$ and $p\to\sigma(p)$ are differentiable on $(p_c,1)$.
\ecor
\begin{proof}
As shown in \cite[(2.4)]{SS13},
\be\label{def-sigma-alpha}
\alpha(p) \ = \ \frac{\E_p[X_2]}{\E_p[\tau_2]} \ \ \mbox{and} \ \sigma^2(p) \ = \frac{\E_p\left[ (X_2\E_p[\tau_2]-\tau_2\E_p[X_2])^2  \right]}{\E_p[\tau_2]^3}.
\ee
The result then follows from a direct application of Proposition \ref{lem:f-g-h-diff} below.
\end{proof}

\begin{proof}[Proof of Proposition \ref{lem:f-g-h-diff}]
Kuczek \cite{K89} showed that when $p>p_c$, the distribution of $\tau_2$ has exponential tail. Since $|X_2|\leq \tau_2$, it follows that the infinite series
\begin{equation}
f_{i,j}(p) \ := \ \sum_{(x,n)\in\Z\times\N: x\in[-n,n]} \ x^i n^j  \ \P_p( X_2=x, \ \tau_2 = n )
\end{equation}
converges for each $p\in (p_c, 1)$. Let $C_o(n)$ be the set of edges (closed or open) discovered before exploring
time horizon $n$ (i.e. before reaching the point $(\rho_o(n),n)$). Again from \cite{K89}, we know that
$$
\P_p( X_2=x, \ \tau_2 = n ) \ = \ \P_p\left(C_o(n)\in A_{x,n}\right) 
$$
where we say $C_o(n)\in A_{x,n}$ if $\rho_o$ satisfies the following property:
\begin{equation}\label{CoAxn}
o\to(\rho_o(n),n), \ \ \rho_o(n)=x, \mbox{   and   }  (\rho_o(i),i)\not\to (x,n)\ \forall\, 1 \leq i<n.
\end{equation}
Note that this property is entirely determined by the realization of $C_o(n)$.

To prove the differentiability of
$$
f_{i,j}(p) = \sum_{(x,n)\in\Z\times\N: x\in[-n,n]} \ x^i n^j  \ \P_p\left(C_o(n)\in A_{x,n}\right),
$$
we will first show that (i) $\P_p\left(C_o(n)\in A_{x,n}\right)$ is differentiable on $(p_c,1]$ for every pair $(x,n)$; and then show that (ii) the series
\begin{equation}\label{derseries}
 \sum_{(x,n)\in\Z\times\N: x\in[-n,n]} \ x^i n^j \frac{\partial }{\partial p} \P_p\left(C_o(n)\in A_{x,n}\right)
\end{equation}
converges uniformly on every interval $[a,b]\subset (p_c,1)$.

First note that
\begin{equation}
\P_p\left(  C_o(n)\in A_{x,n} \right) \ = \ \sum_{C_o(n)\in A_{x,n}} \ p^{|C^{\rm open}_o(n)|} (1-p)^{|C^{\rm closed}_o(n)|},
\end{equation}
where $C^{\rm open}_o(n)$ and $C^{\rm closed}_o(n)$ denote respectively the subset of open and closed edges in $C_o(n)$. Note that the sum above is finite since under the condition $C_o(n)\in A_{x,n}$, we have
$$
C_o(n)\subseteq\{ e = \left((x,t),(x\pm1,t+1)\right) \ : \ 0\leq t\leq n, \ |x|\leq n\}.
$$
It follows that $\P_p\left(C_o(n)\in A_{x,n}\right)$ is differentiable with
\begin{equation}
\frac{\partial }{\partial p} \P_p\left(C_o(n)\in A_{x,n}\right) = \sum_{C_o(n)\in A_{x,n}} \left( \frac{|C^{\rm open}_o(n)|} {p} \ - \ \frac{|C^{\rm closed}_o(n)|} {1-p}    \right) p^{|C^{\rm open}_o(n)|} (1-p)^{|C^{\rm closed}_o(n)|},
\end{equation}
which implies that for all $p\in [a,b]$,
\begin{equation}
\begin{aligned}
\Big| \frac{\partial }{\partial p} \P_p\left(C_o(n)\in A_{x,n}\right) \Big| & \leq \Big(\frac{1}{a}+\frac{1}{1-b}\Big) \sum_{C_o(n)\in A_{x,n}} |C_o(n)| p^{|C^{\rm open}_o(n)|} (1-p)^{|C^{\rm closed}_o(n)|} \\
& \leq c n^2 \P_p(C_o(n)\in A_{x,n}) \leq  c n^2  \P_p(\tau_2 = n),
\end{aligned}
\end{equation}
where we used $|C_o(n)|\leq 2n(2n+1)$, and the constant $c$ depends only on $a$ and $b$. It follows that
$$
 \sum_{(x,n)\in\Z\times\N: x\in[-n,n]} \ |x|^i n^j \Big|\frac{\partial }{\partial p} \P_p\left(C_o(n)\in A_{x,n}\right)\Big|
\leq c \cdot 3^i \sum_{n\in\N} n^{i+j+3} \P_p(\tau_2=n).
$$
To prove the uniform convergence of the series in \eqref{derseries} for $p\in [a,b]$, we only need to prove the uniform convergence of
\begin{equation}
\sum_{n\in\N} n^{i+j+3} \P_p(\tau_2=n), \qquad p\in [a,b].
\end{equation}
For each $p>p_c$, Kuczek~\cite{K89} has proved that $\P(\tau_2\geq n) \leq C_1(p) e^{-C_2(p)n}$, and hence the above series converges pointwise.
The uniform convergence follows from the fact that $C_1(p)$ and $C_2(p)$ can be chosen uniformly for $p\in [a, b]$, which also follows from Kuczek's proof. The key ingredient is the following estimate by Durrett~\cite{D84}:
\begin{equation}\label{decayest}
\P_p( o\to \Z\times \{n\}, o\notin {\cal K}) \leq c_1 e^{-c_2 n},
\end{equation}
where it is easily seen from the proof that $c_1, c_2$ can be chosen uniformly for $p\in [a,b]\subset (p_c, 1)$. This concludes the proof of Proposition \ref{lem:f-g-h-diff}.
\end{proof}

\section{Invariance principle for a pair of left-right paths}\label{sect:single-pair}
Following the notation in \cite{SS13}, let $\Pi$ denote the space of real-valued continuous functions starting from some time in $\R$, equipped with the topology of local uniform convergence plus convergence of the starting time.
Recall from Section \ref{intro} the definition of the percolation configurations $\Omega^{\eps,\pm}$, and for $z\in \Z^2_{\rm even}$, let $\gamma^{\epm}_z$ and $\rho^{\epm}_z$ be paths defined from $\Omega^{\eps,\pm}$ in the same way as $\gamma_z$ and $\rho_z$ are defined from $\Omega$ in Section \ref{sect:preliminaries}.
The main result of this section is the following.

\bteo[Convergence of a pair of left-right paths]\label{teo:cond-conv-sticky-BM}
Let $z_\pm:=(x_\pm, t_\pm)\in \R^2$ and let $z_\pm^\eps:=(x^\eps_\pm, t^\eps_\pm) \in \Z^2_{\rm even}$ be such that $S_{\alpha(p), \sigma(p), \eps} (z^\eps_\pm) \to z_\pm$ as $\eps\to 0$. Let $(l,r)$ be a left-right pair of sticky Brownian motions with drift $b=\frac{\alpha'(p)}{\sigma(p)}$, starting respectively at the space-time points $z_-$ and $z_+$. Then
\begin{equation}\label{pairlimit}
\Big(S_{\alpha(p),\sigma(p),\eps}(\rho^{\eps,-}_{z^\eps_-}, \gamma^{\eps, -}_{z^\eps_-}, \rho^{\eps,+}_{z^\eps_+}, \gamma^{\eps, +}_{z^\eps_+})\Big) \Astoo{\eps} \left(l, l, r, r\right),
\end{equation}
where $\rho^{\eps, \pm}$ are defined at non-integer times via linear interpolation, and the scaling map $S_{\alpha(p),\sigma(p),\eps}$ is defined in \eqref{Seps}.
\eteo

In the next section, we provide an outline the proof. We postpone the main technical parts until further sections.

\subsection{Outline of the proof}
Until further notice, we will assume that $t^\eps_-=t^{\eps}_+$. In the final subsection \ref{sect:diff-starting-times},
we will show how Theorem \ref{teo:cond-conv-sticky-BM} can easily be deduced from this particular case.

We decompose the proof into several steps.
Recall from the introduction that the pair of left-right sticky pair $(l, r)$ with drift $b$ is defined as the unique weak solution of the {\em left-right SDE}
\begin{equation}\label{lrsde}
\begin{aligned}
d l(t) & = 1_{l_t\neq r_t} dB^l(t) + 1_{l_t=r_t} dB^s(t) - bdt \\
d r(t) & = 1_{l_t\neq r_t} dB^r(t) + 1_{l_t=r_t} dB^s(t) + bdt
\end{aligned}
\end{equation}
subject to the constraint that $l(t)\leq r(t)$ for all $t\geq \tau:=\inf\{s: l(s)\leq r(s)\}$, and $B^l, B^r, B^s$ are independent standard Brownian motions. We first need the following characterization of $(l,r)$.

\bprop\label{prop:lr}
Let $b>0$, and let $z_\pm=(x_\pm, t_\pm)$. The sticky pair $(l,r)$ is the unique process satisfying the following three properties:
\begin{enumerate}
\item[\rm (a)] $l=(l(t))_{t\geq t_-}$ and $r=(r(t))_{t\geq t_+}$ are Brownian motions (defined w.r.t.\ the same filtration) with unit diffusion constant, respective drift $-b$ and $b$, and $l(t_-)=x_-$, $r(t_+)=x_+$.
\item[\rm (b)] $1_{l(t) \neq r(t)}\, d\langle l,r\rangle_t = 0$, where $\langle l, r\rangle$ denotes the cross-variation of $l$ and $r$.
\item[\rm (c)] $l(t) \leq r(t)$ for all $t\geq \tau:=\inf\{s\geq t_-\vee t_+: l(s)\leq r(s)\}$.
\end{enumerate}
Furthermore, $\frac{1}{\sqrt{2}}(r(\tau+t)-l(\tau+t))$ is equal in law to the weak solution of the following SDE:
\be\label{eq:sticky-bm}
dw(t) \ = \ \sqrt{2} b \, dt + 1_{w\neq 0}\, dB(t), \qquad w(0)=\frac{r(\tau)-l(\tau)}{\sqrt 2},
\ee
where $B$ is a standard Brownian motion.
\eprop
\noindent
We postpone the proof till Section \ref{sect:proof-lr}.
\bigskip

We first establish the invariance principle for $\gamma$ and $\rho$ starting from a single point, for which it suffices to consider paths $\rho^{\eps, \pm}:= \rho_o^{\eps,\pm}$ starting at the origin. Denote by $(\rho^{\eps,\pm}(T_i^{\eps,\pm}),T_i^{\eps,\pm})_{i\in\N}$ the successive break points along $\rho^{\eps,\pm}$, and let $(X_i^\epm,\tau_i^\epm)_{i\in\N}$ be the successive space-time increments between break points as defined in Section \ref{sect:preliminaries}.

\blem\label{lem:W-T}
For every $n\in\N$, define
\beqn
W^\epm(n) \ := \ \sum_{i=1}^n  \left(X_i^\epm - \alpha(p) \tau_i^\epm \right),  \  \ T^\epm(n) \ := \ T_{n}^{\epm} =  \sum_{i=1}^n \tau^\epm_i \label{def:W-T},
\eeqn
and for $t\notin\N$, define $W(t),T(t)$ by linear interpolation. Conditional on $o\in {\cal K}^{\epm}$,
\begin{align}
\left(\frac{\eps}{\sigma(p)} W^{\epm}(\eps^{-2}t)\right)_{t\geq 0}\quad  & \Astoo{\eps}\quad  \left(B(\E[\tau_2]t) \pm \E[\tau_2] \frac{\alpha'(p)}{\sigma(p)} t\right)_{t\geq0},\label{110} \\
(\eps^2T^\epm(\eps^{-2}t))_{t\geq 0} \quad & \Astoo{\eps}\quad \left(t \E[\tau_2]\right)_{t\geq0}, \label{111}
\end{align}
where $B$ is a standard Brownian motion.
\elem

\begin{proof}
From \cite{K89}, conditioned on $o\in{\cal K}^\epm$, $(X_i^\epm,\tau_i^\epm)_{i\in\N}$ is a sequence of i.i.d random variables with the same distribution
as $(X_2^\epm, \tau_2^\epm)$ without conditioning. Furthermore, using Proposition \ref{lem:f-g-h-diff} we have
\begin{equation}
\begin{aligned}
\E\left[X_2^\epm - \alpha(p) \tau_2^\epm\right] \ &= \ \pm \ \E[\tau_2] \alpha'(p)\eps + O(\eps^2), \\
 \ \ \E\left[(X_2^\epm - \alpha(p) \tau_2^\epm)^2\right] \ &= \ \E[\tau_2] \sigma^2(p) + O(\eps).
\end{aligned} \label{eq:approx}
\end{equation}
The first convergence statement then follows from a standard invariance principle (for triangular arrays).
The second convergence statement follows from the law of large numbers (for triangular arrays) and the fact that
$T^\epm$ is non-decreasing.
\end{proof}

\bcor[Conditional Invariance Principle]\label{cor:conv-marginal}

Conditional on $o\in {\cal K}^{\eps,-}$, resp., ${\cal K}^{\eps,+}$,
$$
S_{\alpha(p),\sigma(p),\eps}(\rho^{\eps,-}) \Astoo{\eps} l, \quad  \mbox{resp.,} \quad S_{\alpha(p),\sigma(p),\eps}(\rho^{\eps,+}) \Astoo{\eps} r,
$$
where $l$, resp., $r$, is distributed as in the left-right pair of sticky Brownian motions with drift $b=\frac{\alpha'(p)}{\sigma(p)}$, defined in Proposition \ref{prop:lr}.
\ecor

\begin{proof}
We follow closely the proof Proposition 2.2 in \cite{SS13}. First, we start by replacing
the path $\rho^\epm$
by $\tilde \rho^\epm$ where $\tilde \rho^\epm(T_i^\epm) = \rho^\epm(T_i^\epm)$ and for $t\in(T_i^\epm,T_{i+1}^\epm)$, $\tilde \rho^\epm$
is defined by linear interpolation.
More precisely, we write the decomposition:
$$
S_{\alpha(p),\sigma(p),\eps}(\rho^\epm) \ = \ S_{\alpha(p),\sigma(p),\eps}(\tilde \rho^\epm) \ + \ S_{0,\sigma(p),\eps}(\tilde \rho^\epm - \rho^\epm).
$$
Arguments analogous to the ones in \cite{SS13} show that the second term vanishes as $\eps\to0$.

From the definition, it is straightforward to see that
$$
S_{\alpha(p),0,0}\left( \tilde \rho^\epm \right)(t) \ = \ W^\epm((T^\epm)^{-1}(t)),
$$
and
$$
S_{\alpha(p),\sigma(p),\eps}(\tilde \rho^\epm) = S_{0,\sigma(p),\eps} S_{\alpha(p),0,0}(\tilde \rho^\epm) = S_{0,\sigma(p),\eps}(W^\epm((T^\epm)^{-1})).
$$
Lemma \ref{lem:W-T} implies that $S_{0,\sigma(p),\eps}(W^\epm)$ converges to
$(B(\E[\tau_2]t) \pm \E[\tau_2] \frac{\alpha'(p)}{\sigma(p)}t)_{t\geq0}$, where $B$
is a standard Brownian motion, and $({\eps^2}T^\epm(t/\eps^2))_{t\geq0}$
converges to $(\E[\tau_2] t)_{t\geq0}$. It then follows that $S_{\alpha(p),\sigma(p),\eps}(\rho^{\epm})$ converges in distribution to
$$
\Big(B(t) \ \pm \ \frac{\alpha'(p)}{\sigma(p)}t\Big)_{t\geq0},
$$
which is exactly the law of $l$, resp., $r$.
\end{proof}

\bprop[Invariance Principle]\label{prop:uncond-inv-principle} Following the same notation as in Corollary \ref{cor:conv-marginal}, we have
$$
S_{\sigma(p),\alpha(p),\eps}(\rho^{\eps,-}, \gamma^{\eps, -}) \Astoo{\eps} (l, l), \quad \mbox{and} \quad S_{\sigma(p),\alpha(p),\eps}(\rho^{\eps,+}, \gamma^{\eps, +}) \Astoo{\eps} (r, r).
$$
\eprop
\begin{proof}
The unconditional invariance principle for $\rho^{\eps,\pm}$, as well as the fact that $\rho^{\eps, \pm}$ and $\gamma^{\eps, \pm}$ converge to the same scaling limit, can be established by the same argument as in the proof of \cite[Prop. 2.2]{SS13}. The only modification needed is that we need to use the fact that the tail estimate \eqref{decayest} is uniform in $p\in [a,b]\subset (p_c, 1)$, since the percolation parameters now depend on $\eps$.
\end{proof}

To establish Theorem \ref{teo:cond-conv-sticky-BM}, note that Proposition \ref{prop:uncond-inv-principle} implies that $(S_{\sigma(p),\alpha(p),\eps}(\rho^{\eps,-}_{z^\eps_-}, \rho^{\eps,+}_{z^\eps_+}))_{\eps>0}$ is a tight sequence of $\Pi\times\Pi$-valued random variables. Going to a subsequence if necessary, we can assume the existence of a pair of drifted Brownian motions $(l,r)$ such that
$$
S_{\sigma(p),\alpha(p),\eps}(\rho^\emi_{z^\eps_-}, \rho^\epl_{z^\eps_+}) \Astoo{\eps} (l,r).
$$
In order to prove Theorem \ref{teo:cond-conv-sticky-BM}, it only remains to show that the sub-sequential limit is unique in distribution, i.e.,  $(l, r)$ satisfies the three conditions of Proposition \ref{prop:lr}. Condition $(a)$ is satisfied by Proposition \ref{prop:uncond-inv-principle}.
Conditions (b) and (c) will follow from the next two results.

\bprop\label{41}
$1_{l(t)\neq r(t)} d \langle l,r\rangle_t = 0$.
\eprop

\bprop\label{cond-1-2}
$l(t)\leq r(t)$ for all $t\geq \tau:=\inf\{s\geq t_-\vee t_+: l(s)\leq r(s)\}$.
\eprop

In order prove Theorem \ref{teo:cond-conv-sticky-BM}, it remains to show  Propositions \ref{prop:lr}, \ref{41}, and \ref{cond-1-2}. This is done in the next sections.


\subsection{Proof of Proposition \ref{prop:lr}}\label{sect:proof-lr}

It has been shown in \cite{SS08} that the left-right SDE \eqref{lrsde} has a unique weak solution $(l, r)$, and $\frac{1}{\sqrt 2}(r-l)$ is a weak solution of the equation \eqref{eq:sticky-bm}.  Any such solution $(l, r)$ clearly satisfies conditions (a)--(c) in Proposition \ref{prop:lr}.

\bigskip
It remains to prove uniqueness of the process satisfying the three properties.
Let us assume that there exists a process $(l,r)$ satisfying Proposition \ref{prop:lr} (a)-(c). Note that (a) and (b) imply that $l$ and $r$ are independent Brownian motions before they meet at time $\tau$. Therefore we only need to prove uniqueness in distribution from time $\tau$ onward, and hence we may assume that $l$ and $r$ both start at $0$ at time $0$. Define
$$
U(t) := \frac{r(t)+l(t)}{2}, \qquad V(t):= \frac{r(t)-l(t)}{2} - bt.
$$
By assumption (a), it is straightforward to check that $U$ and $V$ are orthogonal martingales. Furthermore, writing
\begin{equation}\label{Ct}
C(t)=\frac{t}{2}- \frac{1}{2}\int_0^t 1_{l(s) = r(s)} d\langle l,r\rangle_s,
\end{equation}
we have $\langle U\rangle(t)=t-C(t)$ and $\langle V\rangle (t)=C(t)$. By Knight's theorem~\cite[Theorem V~(1.9)]{RY99}, there exist two independent standard Brownian motions $B$ and $B'$ such that
$$
U = B'(t-C(t)), \qquad V(t) =  B(C(t)).
$$
Recovering $l$ and $r$ from $U$ and $V$, we obtain
\begin{equation}\label{lrCt}
\begin{aligned}
r(t) & =  B'(t-C(t)) + B(C(t)) + b t, \\
l(t) & =  B'(t-C(t)) - B(C(t)) - b t.
\end{aligned}
\end{equation}

Note that $C$ is increasing. We claim that $C$ is strictly increasing. Assume the contrary that $C(t_1)=C(t_2)$ for some $0\leq t_1<t_2$. Then
by \eqref{Ct}, we have
$$
\frac{t_2-t_1}{2} = \frac{1}{2}\int_{t_1}^{t_2} 1_{l(s) = r(s)} d\langle l,r\rangle_s,
$$
which is only possible if $l(s)=r(s)$ for all $s\in [t_1, t_2]$, since $d\langle l, r\rangle_s\leq ds$. However, by \eqref{lrCt}, we cannot
have both $l=r$ and $C$ being a constant on $[t_1, t_2]$. Therefore $C(t)$ is strictly increasing and admits an inverse $C^{-1}$. Denote
$$
Z(t) := \frac{r(C^{-1}(t))- l(C^{-1}(t))}{2} = B(t) + bC^{-1}(t).
$$
In order to show the uniqueness the law of $(l, r)$, we will
show:
\begin{itemize}
\item[(1)] $Z$ is $B(t)+2bt$ Skorohod reflected at $0$;
\item[(2)] $C^{-1}(t)=2t - \frac{1}{b} \inf_{s\in [0,t]} (B(s)+2bs)\wedge 0$.
\end{itemize}
This would imply that $(C(t))_{t\geq 0}$ is determined by $(B(t))_{t\geq 0}$, and hence the law of $(l, r)$ is unique by \eqref{lrCt}.
We can write
$$
Z(t) = B(t) + 2bt + b (C^{-1}(t) - 2t).
$$
By \eqref{Ct}, it is easily seen that $(C^{-1}(t) - 2t)$ is non-negative, and it only increases when $Z=0$. Since $Z\geq0$ by Proposition \ref{prop:lr}~(c), by uniqueness of the solution to the Skorohod equation~\cite[Lemma VI~(2.1)]{RY99}, we deduce that
$Z$ is the Skorohod reflection of $B(t)+2bt$ at $0$, and
$$
b (C^{-1}(t) - 2t) = -\inf_{s\in [0,t]} (B(s)+2bs)\wedge 0.
$$
This completes the proof of (1) and (2), and thus the proof of uniqueness.
\qed

\subsection{Proof of Proposition \ref{41}}\label{sect:proof-41}

The following preliminary result will be needed to prove Proposition \ref{41}. It shows that two paths $\rho^{\eps, \pm}_{z_\pm}$ converge in the scaling limit to two independent Brownian motions before they meet.

\blem\label{prop:conv-pair-expl-cluster}
Let $z_\pm:=(x_\pm, t_\pm)\in \R^2$ and let $z_\pm^\eps:=(x^\eps_\pm, t^\eps_\pm) \in \Z^2_{\rm even}$ be such that $S_{\alpha(p), \sigma(p), \eps} (z^\eps_\pm) \to z_\pm$ as $\eps\to 0$. Assume that $z_+\neq z_-$, and for any $\delta>0$, define
$$
\tau^\eps_\delta = \inf\Big\{n\geq t^\eps_+ \vee t^\eps_- \ : \ \frac{\eps}{\sigma(p)} \left|\rho^{\epl}_{z^\eps_+}(n) - \rho^\emi_{z^\eps_-}(n)\right| \ \leq \ \delta \Big\}.
$$
Then
\begin{align*}
&\left(S_{\alpha(p),\sigma(p),\eps}\left(  \rho^\emi_{z^\eps_-}(\cdot\wedge \tau_\delta^\eps), \, \gamma^\emi_{z^\eps_-}(\cdot\wedge \tau_\delta^\eps), \, \rho^{\epl}_{_{z^\eps_+}}(\cdot\wedge \tau_\delta^\eps), \, \gamma^{\epl}_{_{z^\eps_+}}(\cdot\wedge \tau_\delta^\eps) \right), \, \eps^2 \tau_\delta^\eps\right) \\
& \hspace{6cm} \Astoo{\eps}  \quad  \left( l(\cdot\wedge \tau_\delta), \, l(\cdot\wedge \tau_\delta), \, r(\cdot\wedge \tau_\delta), \, r(\cdot\wedge \tau_\delta),  \, \tau_\delta \right),
\end{align*}
where $(l,r)$ is a pair of independent Brownian motions starting respectively at $z_-$ and $z_+$, with respective drift $-b$ and $b$ and common diffusion constant $1$, and
$$
\tau_\delta := \inf \{t\geq t_+\vee t_- : |r(t)-l(t)|< \delta \}.
$$
\elem

\begin{proof} The proof is the same as in the proof of \cite[Prop. 3.3]{SS13}. The main idea is that the paths $\rho^{\eps, \pm}_{z^\eps_\pm}$ and $\gamma^{\eps, \pm}_{z^\eps_\pm}$ can be approximated by their respective percolation exploration clusters, which evolve independently before they intersect. Since the exploration clusters (which are always bounded between $\rho_z$ and $\gamma_z$) individually converge to a Brownian motion by Proposition \ref{prop:uncond-inv-principle}, the result then follows. See the proof of \cite[Prop. 3.3]{SS13} for more details, which is more involved than our case since it also considers what happens after the exploration clusters intersect.
\end{proof}

By going to a subsequence if needed, we may assume that $S_{\sigma(p),\alpha(p),\eps}(\rho^\emi_{z^\eps_-}, \rho^\epl_{z^\eps_+})\Rightarrow (l, r)$. To prove $1_{l(t)\neq r(t)} d\langle l, r\rangle_t=0$, it suffices to show that for any $\delta>0$, $(l, r)$ are distributed as independent Brownian motions when $\Delta(t):= r(t)-l(t)$ satisfies $|\Delta(t)|>\delta$.

Define
$
\kappa_0,\xi_0=t_-\vee t_+,
$
and for every $i\geq0$:
\begin{equation}\label{kappaxi}
\kappa_{i+1}=\inf\{t\geq\xi_i \ : |\Delta(t)| > \delta\}, \qquad
\xi_{i+1}=\inf\{t\geq\kappa_{i+1} \ : \ |\Delta(t)| < \delta/2 \},
\end{equation}
with the convention that $\inf\{\emptyset\}=\infty$. To prove Proposition \ref{41}, it then suffices to show that
\blem\label{lem:reg-and-BM}
For each $i\in\N$, conditioned on $(l(t), r(t))_{t\leq \kappa_i}$ with $\kappa_i<\infty$, the process
$$
\big(l(s\wedge \xi_i), \, r(s\wedge \xi_i)\big)_{s\geq \kappa_i}
$$
is distributed as a pair of independent Brownian motions with unit diffusion constant and respective drift $-b, b$, stopped when their distance reaches $\delta/2$.
\elem
\begin{proof}
This will be proved by an approximation argument. Recall that we have assumed that
$$
S_{\sigma(p),\alpha(p),\eps}(\rho^\emi_{z^\eps_-}, \rho^\epl_{z^\eps_+}) \ \Astoo{\eps} \ (l, r).
$$
Analogous to $(\kappa_i,\xi_i)$, let us define $\kappa^\eps_0=\xi^\eps_0 = t^\eps_-\vee t^\eps_+$, and for $i\geq1$, define
$$
 {\kappa}_{i+1}^\eps = \inf\Big\{ k > \xi_i^\eps \ : \ |\Delta \rho^{\eps}(k)| \geq \frac{\sigma(p)}{\eps}\delta \Big\}, \qquad
 {\xi}_{i+1}^\eps = \inf\Big\{ k > \kappa_{i+1}^\eps \ : \ |\Delta \rho^{\eps}(k)| \leq \frac{\sigma(p)}{\eps}\frac \delta 2  \Big\},
$$
where $\Delta \rho^{\eps} := \rho^\epl_{z^\eps_+}-\rho^\emi_{z^\eps_-}$. As $[0, \infty]$-valued random variables, $(\eps^2{\kappa}_i^\eps)_{i\in\N}$ and $(\eps^2{\xi}_i^\eps)_{i\in\N}$ are automatically tight. Therefore by going to a subsequence if necessary, we may assume that
\be\label{eq:conv-succ-times}
\left(S_{\alpha(p),\sigma(p),\eps}\left(\rho^\emi_{z^\eps_-}, \gamma^\emi_{z^\eps_-}, \rho^\epl_{z^\eps_+}, \gamma^\epl_{z^\eps_+}\right), \ (\eps^2{\kappa}_i^\eps,  \eps^2{\xi}_i^\eps)_{i\in\N} \right) \Longrightarrow \left((l, l ,r, r), \ (\bar \kappa_i, \bar \xi_i)_{i\in\N} \right),
\ee
and by the Skorohod representation theorem, we can assume the above convergence to be almost sure via a suitable coupling. Note that under this coupling, by the definition of $(\kappa_i, \xi_i)_{i\in\N}$ in \eqref{kappaxi}, we must have $\bar\kappa_i \leq \kappa_i$ and $\bar\xi_i\leq \xi_i$ for all $i\in\N$, which will be strengthened to equalities later.

Let us fix $i\in\N$, and to simplify notation, we will omit the subscript $i$ from $\kappa_i^\eps, \bar \kappa_i, \xi_i^\eps$. Let us consider the composite path $\bar \rho^\epm$ defined as follows:
\begin{equation}\label{324}
\bar \rho^\epm_{z^\eps_\pm}(t) \ = \ \left\{
\begin{aligned}
\rho^\epm_{z^\eps_\pm}(t) &\qquad  \mbox{if $t\leq \kappa^\eps$,} \\
\rho_{w^\eps_\pm}^\epm(t) &\qquad  \mbox{if $t\geq \kappa^\eps$, where $w^\eps_\pm := (  \rho_{z^\eps_\pm}^\epm(\kappa^\eps), \ \kappa^\eps)$.}
\end{aligned} \right.
\end{equation}
We first claim that conditional on $\{\bar \kappa<\infty\}$, these modified paths are good approximations of $\rho^\epm_{z^\eps_\pm}$, and hence of $(l, r)$. More precisely,
\be\label{eq:approx2-(r,l)}
S_{\alpha(p),\sigma(p),\eps}\big(\rho^{\eps, \pm}_{w^\eps_\pm}(t) - \rho_{z^\eps_\pm}^{\eps, \pm}(t)\big)_{t\geq \kappa^\eps} \ \Astoo{\eps} \ 0.
\ee
Note that conditioned on $\bar \kappa<\infty$, $S_{\alpha(p),\sigma(p),\eps}(\rho^{\eps, \pm}_{w^\eps_\pm}- \gamma^{\eps, \pm}_{w^\eps_\pm})  \Rightarrow 0$ by Proposition \ref{prop:uncond-inv-principle}. To prove \eqref{eq:approx2-(r,l)}, it then suffices to show that
\begin{equation}\label{326}
S_{\alpha(p),\sigma(p),\eps}\big(\gamma^{\eps, \pm}_{w^\eps_\pm}(t) - \gamma_{z^\eps_\pm}^{\eps, \pm}(t)\big)_{t\geq \kappa^\eps} \ \Astoo{\eps} \ 0.
\end{equation}
This is easily seen to hold since the difference between the starting points of these two paths,
$$
|\rho_{z^\eps_\pm}^\epm(\kappa^\eps)-\gamma_{z^\eps_\pm}^\epm(\kappa^\eps)| \astoo{\eps}{0} 0
$$
almost surely on the event $\{\bar\kappa<\infty\}$ by the coupling assumed in \eqref{eq:conv-succ-times}. In particular, for any $\alpha>0$,
if we denote
$$
\widetilde w^\eps_\pm := (\rho_{z^\eps_\pm}^\epm(\kappa^\eps) - \sigma\alpha/\eps, \ \kappa^\eps),
$$
then except on a set with probability tending to $0$, we have
\begin{equation}\label{327}
\gamma_{\widetilde w^\eps_\pm}^\epm(t) \leq \gamma_{z^\eps_\pm}^\epm(t) \leq \gamma_{w^\eps_\pm}^\epm(t) \qquad \mbox{for all } t\geq \kappa^\eps.
\end{equation}
The advantage of working with $\gamma_{\widetilde w^\eps_\pm}^\epm$ and $\gamma_{w^\eps_\pm}^\epm$ is that except for their starting points, they depend only on the percolation configuration above time $\kappa^\eps$. Therefore by Proposition 3.3 of \cite{SS13} (modified to take into account the $\eps$-dependence of the percolation parameter), $S_{\alpha(p),\sigma(p),\eps}(\gamma_{\widetilde w^\eps_\pm}^\epm, \gamma_{w^\eps_\pm}^\epm)$
can be approximated by two coalescing Brownian motions started at the same time with distance $\alpha$ apart. By \eqref{327} with $\alpha>0$ chosen arbitrarily small, we then obtain \eqref{326}, and hence \eqref{eq:approx2-(r,l)}.

By going to a further subsequence if necessary, we can now enhance \eqref{eq:conv-succ-times} by including the convergence in \eqref{eq:approx2-(r,l)}, and by coupling via Skorohod representation, we may assume the convergence to be almost sure, so that $S_{\alpha(p), \sigma(p), \eps}(\rho^{\eps, -}_{z^\eps_-}(t), \rho^{\eps, +}_{z^\eps_+}(t))_{t\leq \kappa^\eps}$ converges almost surely to $(l(t), r(t))_{t\leq \bar\kappa}$, and $S_{\alpha(p), \sigma(p), \eps}(\rho^{\eps, -}_{w^\eps_-}, \rho^{\eps, +}_{w^\eps_+})$ converges almost surely to $(l(t), r(t))_{t\geq \bar\kappa}$. On the other hand, we note that conditioned on the percolation configuration up to time $\kappa^\eps$, the law of $(\rho^{\eps, -}_{w^\eps_-}, \rho^{\eps, +}_{w^\eps_+})$ depends only on their starting points $w^\eps_\pm$. Therefore under the coupling that
$S_{\alpha(p), \sigma(p), \eps}(\rho^{\eps, -}_{z^\eps_-}(t), \rho^{\eps, +}_{z^\eps_+}(t))_{t\leq \kappa^\eps}$ converges almost surely to $(l(t), r(t))_{t\leq \bar\kappa}$, by Lemma \ref{prop:conv-pair-expl-cluster}, the law of $S_{\alpha(p), \sigma(p), \eps}(\rho^{\eps, -}_{w^\eps_-}(t\wedge \xi_\eps), \rho^{\eps, +}_{w^\eps_+}(t\wedge \xi^\eps))_{t\geq \kappa^\eps}$ conditioned on $(\rho^{\eps, -}_{z^\eps_-}(t), \rho^{\eps, +}_{z^\eps_+}(t))_{t\leq \kappa^\eps}$ converges to the law of a pair of independent Brownian motions $(\tilde l, \tilde r)$ with unit diffusion constant and respective drift $-b$ and $b$, starting respectively at $l(\bar\kappa)$ and $r(\bar\kappa)$ at time $\bar\kappa$, and stopped when their distance reaches $\delta/2$. Since conditioned on $(l(t), r(t))_{t\leq \bar\kappa}$, $(l(t\wedge \bar\xi), r(t\wedge \bar\xi))_{t\geq \bar\kappa}$ has the same distribution as $(\tilde l, \tilde r)$, Lemma \ref{lem:reg-and-BM} would follow once we show that $(\bar\kappa_i, \bar\xi_i)=(\kappa_i, \xi_i)$.

Indeed, Lemma \ref{prop:conv-pair-expl-cluster} implies that for each $i\in\N$, $\bar \xi_i = \inf\{t\geq \bar \kappa_i : |r(t)-l(t)| < \frac \delta 2\}$. Therefore it only remains to show that $\bar \kappa_i=\kappa_i$ for every $i\in\N$. We can proceed inductively. We may assume that
$|r(\xi_0)-l(\xi_0)|<\delta$. By our definition of $\kappa_1$ and the almost sure coupling in \eqref{eq:conv-succ-times}, we must have $\bar\kappa_1\leq \kappa_1$. On the other hand, we have shown that $(l(t\wedge \bar\xi_1), r(t\wedge\bar\xi_1))_{t\geq \bar\kappa_1}$ is distributed as a pair of independent Brownian motions with initial distance $\delta$ at time $\bar\kappa_1$. Therefore $\bar\kappa_1=\inf\{t>\bar\kappa_1: |r(t)-l(t)|>\delta\}$, which implies that $\kappa_1\leq \bar\kappa_1$. Therefore $\bar\kappa_1=\kappa_1$. We can now iterate the argument to show that $\bar\kappa_i=\kappa_i$ for all $i\in\N$.
\end{proof}

\subsection{Proof of Proposition \ref{cond-1-2}}\label{Sec:cond-1-2}
By going to a subsequence if needed, we may assume that $S_{\sigma(p),\alpha(p),\eps}(\rho^\emi_{z^\eps_-}, \rho^\epl_{z^\eps_+})\Rightarrow (l, r)$, and by the Skorohod representation theorem, we may even assume the convergence to be almost sure.

If $\rho^\emi_{z^\eps_-}(t^\eps)\leq \rho^\epl_{z^\eps_+}(t^\eps)$, then we have $\rho^{\eps, -}_{z^\eps_-}\leq \rho^{\eps, +}_{z^\eps_+}$ for all time since $\rho^{\eps, \pm}_{z^\eps_\pm}$ are defined from the percolation configurations $\Omega^{\eps, \pm}$, and almost surely $\Omega^{\eps, +}$ contains more open edges than $\Omega^{\eps, -}$. It follows that $l\leq r$ almost surely.

If $\rho^\epl_{z^\eps_+}(t^\eps) < \rho^\emi_{z^\eps_-}(t^\eps)$, then from the definition of $\rho^{\eps, \pm}_{z^\eps_\pm}$, it is easy to see that $\rho^\epl_{z^\eps_+}(n) \geq  \rho^\emi_{z^\eps_-}(n)$ for all $n$ larger than or equal to
$$
\tau^\eps_\rho:=\inf\{i\geq t^\eps: \rho^\epl_{z^\eps_+}(i) \geq  \rho^\emi_{z^\eps_-}(i) \}.
$$
Following the same arguments as in the proof of \cite[Lemma 3.1]{SS13}, we note that $\tau^\eps_\rho$ is the first time the percolation cluster $C_{z^\eps_+}$ in the percolation configuration $\Omega^{\eps, +}$ intersects the percolation cluster $C_{z^\eps_-}$ in the percolation configuration $\Omega^{\eps, -}$. Since the percolation clusters evolve independently before they intersect, and each converges to a Brownian motion with respective drift $\pm b$, it is easily seen (cf. proof of \cite[Prop.~3.3]{SS13}) that $\eps^2 \tau^\eps_\rho$ converges in distribution to $\tau:=\inf\{s\geq t_-\vee t_+: l(s)\leq r(s)\}$. Since $\rho^\epl_{z^\eps_+}(n) \geq  \rho^\emi_{z^\eps_-}(n)$ for $n\geq  \tau^\eps_\rho$, it follows that $r(t)\geq l(t)$ for all $t\geq \tau$, which concludes the proof of Proposition \ref{cond-1-2}.
\qed

\subsection{Different starting times}\label{sect:diff-starting-times}
Until now, we only considered the case $t^\eps_-= t^\eps_+$. Let us now show how to extend this result to the case $t^\eps_-< t^\eps_+$. (The case $t^\eps_- > t^\eps_+$. can be treated along the same lines.) The proof goes by an approximation argument analogous to the proof Lemma \ref{lem:reg-and-BM}. It is sufficient to decompose the path starting from earlier time into two time-windows: before and after reaching the starting time of the later path.

More precisely, we can approximate the path $l(\cdot \wedge t_+)$ by
the path $S_{\alpha(p),\sigma(p),\eps}(\rho^{-,\eps}_{z^\eps_-})$ stopped at time $t^\eps_+$, and the remaining section of the path
(i.e., the path $(l(t); t\geq t_{+}^\eps))$) by the rescaled exploration cluster
starting from the point $\rho^{-,\eps}_{z^\eps_-}(t_+^\eps)$ at time $t_+^\eps$, which by definition, only uses the percolation configuration above time $t^\eps_+$. Since the percolation configurations before and after time horizon $t^\eps_+$ are independent,
and since the latter two approximations only depend on the lower and upper part of the percolation cluster at time $t^\eps_+$ respectively, it follows that conditioned on $l(t_+)$,
the pair of paths $(l(t),r(t); t\geq t^+)$ is distributed as a pair of sticky paths starting from
$(l({t^+}),r({t^+}))$ so that $(l,r)$ is distributed as a left-right pair of sticky Brownian paths.

\section{Convergence to the left-right sticky Brownian webs}
\label{sticky-webs-proof}

In this section, we prove Theorem \ref{thm:main}. Before going into the proof, we first need to clarify the definition of the filtration $({\mathcal G}_t; t
\in\R)$ alluded to in Theorem  \ref{sticky-webs} . Recall that $\Pi$ denotes the set of paths with a starting point equipped with a norm $d_{\Pi}$
inducing  the topology of local uniform convergence plus convergence of the starting time. See e.g.\ the review article \cite{SSS17} for more details.

For any path with starting point  with time coordinate less than $t$, define the killing operator $K_t$ such that $K_t h(s) \ = \ h(s\wedge t)$.
Define $\theta_t$ the time shift operator such that $\theta_t\circ h(s) \ = \ h(s+t)$.

Let us now consider a Brownian web $\W$ and ${\cal D}$ a dense subset of $\R^2$. Define the $\sigma$-field
$$
{\cal F}_t \ = \ \sigma\big( K_t \circ w_{x_i, t_i} : (x_i, t_i)\in {\cal D}, \ t_i\leq t, w_{x_i, t_i} \in {\cal W} \big),
$$
where $w_{x,t}$ is the (a.s.\ unique) path in ${\cal W}$ starting from the point $(x,t)$.
By a standard path coupling argument, see e.g. \cite[Lemma 3.4]{SS08}, it can be proved that ${\cal F}_t$ is independent of the choice of the set ${\cal D}$. We call $({\cal F}_t; t\geq0)$ the natural filtration associated to $\W$. The common filtration ${\cal G}_t$ referred to in Theorem \ref{sticky-webs} is defined as follows.

\begin{defin}
\label{def:common-filtration}
Let $(\bar \Wl, \bar \Wr)$ be two (possibly drifted) coupled Brownian webs and let
$(\bar {\cal F}^l_t; t\geq0)$ and  $(\bar {\cal F}^r_t; t\geq0)$ be  the natural filtration of $\bar \Wl$ and $\bar \Wr$.
We define ${\cal G}_t \ := \ \bar {\cal F}^l_t \vee \bar {\cal F}^r_t$ as the common filtration of the pair $(\bar \Wl, \bar \Wr)$.
\end{defin}


\bigskip

\noindent

\begin{proof}[Proof of Theorem \ref{thm:main}.] Let us now consider a supercritical oriented percolation model with $p>p_c$. In \cite{SS13}, it was shown that
$$
S_{\alpha(p),\sigma(p),\eps}(\Gamma) \Astoo{\eps} {\cal W},
$$
where $ \Gamma=\{\gamma_z \ : \  z\in{\cal K} \}$ is the set of right-most infinite open paths in the percolation configuration. In previous sections (see Theorem \ref{teo:cond-conv-sticky-BM}), we showed that for every $z^\eps$ such that
$S_{\alpha(p),\sigma(p),\eps}(z^\eps)$ converges to $z$, then $S_{\alpha(p),\sigma(p),\eps}(\gamma^{\pm,\eps}_{z^\eps})$ converges to a drifted Brownian motion with drift $\pm b(p)$. Let us denote by $\Gamma^{\pm,\eps}$
the set of right-most infinite paths in the percolation configuration $\Omega^{\pm,\eps}$. Following the exact same steps as in \cite{SS13}, it can be shown that
\[ S_{\alpha(p),\sigma(p),\eps}(\Gamma^{+,\eps}) \Astoo{\eps} {\cal W}_r,  \ \
S_{\alpha(p),\sigma(p),\eps}(\Gamma^{-,\eps}) \Astoo{\eps} {\cal W}_l, \]
where $\W_r$ (resp., $\W_l$) is a drifted Brownian web with drift $b(p)$ (resp., $-b(p)$). (The latter two convergence statements can be seen as a ``triangular" extension of the results proved in \cite{SS13}.) In particular, this implies that  the sequence of random variables $\{ S_{\alpha(p),\sigma(p),\eps}(\Gamma^{+,\eps}, \Gamma^{-,\eps})  \}_{\eps>0}$ is tight and in order to prove Theorem \ref{thm:main}, it remains to prove that any convergent subsequence must be a sticky pair of left-right Brownian webs. According to Theorem \ref{teo:cond-conv-sticky-BM}, it only remains to prove the third condition of Theorem \ref{sticky-webs}, namely that the two Brownian webs are co-adapted.
This amounts to proving that for any deterministic  $z_1,z_2\in\R^2$, the processes
$(r_{z_1},r_{z_2}), (l_{z_1}, l_{z_2})$ and $(r_{z_1},l_{z_2})$
are Markov processes with respect to the common filtration $({\mathcal G}_t; t\geq0)$ induced by $(\Wl, \Wr)$.
We only show the Markov property $(r_{z_1},r_{z_2})$. The two other processes can be handled by an analogous argument.
%
%
%

%

Let $m\in\N\cup\{\infty\}$ and $f$ and $g$ be continuous functions on
$\Pi^{\otimes 2}$ and $\Pi^{\otimes(2+m)}$ respectively. Let $t
\geq0$ and  let $\{z_i\}_{i=1}^2$
and $\{\bar z_i\}_{i=1}^m$ be two collections of points in $\R^2$ with time coordinates less than $t$.
According to the definition of ${\cal G}_t$, it is enough to prove that
conditioned on $(r_{z_1}(t), r_{z_2}(t))$, the law of $(r_{z_1}(s), r_{z_2}(s))_{s\geq t}$ does not depend on
$(r_{z_1}(s), r_{z_2}(s), (l_{\bar z_j}(s), r_{\bar z_j}(s))_{1\leq j\leq m})_{s<t}$ for any finite $m$
This can be proved by an approximation argument analogous to 
that in Section \ref{sect:diff-starting-times}. Namely, we can approximate each path in $((l_{\bar z_j}(s), r_{\bar z_j}(s))_{1\leq j\leq m})_{s<t}$ by their rescaled exploration clusters up to time $t$, and the same for $(r_{z_1}(s), r_{z_2}(s))_{s< t}$, while we approximate $(r_{z_1}(s))_{s\geq t}$ and $(r_{z_2}(s))_{s\geq t}$ by new exploration clusters that start at time $t$. The desired independence then follows from the independence between percolation configurations above and below time horizon $t$ and passing to the limit as $\eps\downarrow 0$.
\end{proof}

\bigskip

\appendix

\section{Alternative constructions of the left-right pair of sticky Brownian webs} \label{S:sticky}

In this section, we give a proof sketch for Theorem \ref{sticky-webs}. We will show that it is equivalent to the characterization of the left-right Brownian webs given in \cite[Theorem 1.5]{SS08}. We will also discuss the characterization of a pair of $\theta$-coupled sticky webs given in \cite[Theorem 4]{HW09}.

\medskip

In Theorem 1.5 of \cite{SS08}, the pair $({\cal W}_l, {\cal W}_r)$ is uniquely characterized
by condition 1 in Theorem \ref{sticky-webs} plus the finite-dimensional distribution, that is,  the joint distribution of any finite collection of paths $(l_{z_1}, \ldots, l_{z_m}; r_{\bar z_1}, \ldots, r_{\bar z_n})$ with deterministic $z_1, \ldots, z_m, \bar z_1, \ldots, \bar z_n\in \R^2$, where $l_{z_i}$ and $r_{\bar z_j}$ are the almost surely unique path in ${\cal W}_l$ and ${\cal W}_r$, starting at $z_i$ and $\bar z_j$, respectively. In \cite{SS08}, the joint law of $(l_{z_1}, \ldots, l_{z_m}; r_{\bar z_1}, \ldots, r_{\bar z_n})$ is characterized by the properties that: (1) paths evolve independently when they are apart; (2) paths of the same type coalesce when they meet; (3) each pair $(l_{z_i}, r_{\bar z_j})$ interact as a left-right pair of sticky Brownian motions as in Definition \ref{def:lrsde} (see \cite{SS08} for more details). Therefore to deduce Theorem 1.3, it suffices to show that conditions 1-3 in Theorem 1.3 imply that $(l_{z_1}, \ldots, l_{z_m}; r_{\bar z_1}, \ldots, r_{\bar z_n})$ satisfy properties (1)-(3) above.

Indeed, the conditions in  Theorem \ref{sticky-webs} imply that, with respect to the common filtration generated by the collection of paths $(l_{z_1}, \ldots, l_{z_m}; r_{\bar z_1}, \ldots, r_{\bar z_n})$, pairwise, paths of the same type coalesce when they meet, and paths of different types interact as left-right sticky Brownian motions. In particular, for any pair of paths $u, v$ in this collection,
\[1_{u\neq v} d\left<u,v\right> = 0,\]
where {\it the cross-variation is taken with respect to the common filtration} generated by the whole collection. This in turn implies that the paths evolve as independent Brownian motions when they are apart. Therefore properties (1)-(3) above hold, and the pair of sticky Brownian webs $({\cal W}_l, {\cal W}_r)$ characterized by Theorem \ref{sticky-webs} is the same as the left-right Brownian web characterized in \cite[Theorem 1.5]{SS08}.

The above sketch, in particular, the determination of the joint law of $(l_{z_1}, \ldots, l_{z_m}; r_{\bar z_1}, \ldots, r_{\bar z_n})$ from conditions 1-3 in Theorem \ref{sticky-webs}, was rigorously implemented by Howitt and Warren in a slightly different context, where they introduced the notion of co-adaptedness. Their argument can be easily adapted to our setting. In \cite[Theorem 4]{HW09}, they formulated an analogue of Theorem \ref{sticky-webs} for a pair of $\theta$-coupled sticky webs $({\cal W}, {\cal W}')$, where ${\cal W}$ and ${\cal W}'$ each is distributed as a standard Brownian web with no drift, and each pair of path $(\pi, \pi')\in ({\cal W}, {\cal W}')$ with deterministic starting points is distributed as a pair of sticky Brownian motions with stickiness parameter $\theta$.  Again, the key is to show that the analogues of conditions 1-3 in Theorem \ref{sticky-webs} for the $\theta$-coupled sticky webs $({\cal W}, {\cal W}')$ uniquely determine the joint distribution of a finite collection of paths $(\pi_1, \ldots, \pi_m; \pi'_1, \ldots, \pi'_n)$ in $({\cal W}, {\cal W}')$. More precisely, given
$(\pi_1, \ldots, \pi_m; \pi'_1, \ldots, \pi'_n)$, if each pair of paths of the same type are distributed as coalescing Brownian motions with respect to the common filtration, and each pair of paths of different types are distributed as $\theta$-coupled sticky Brownian motions with respect to the common filtration, then the joint distribution of the whole collection is uniquely determined. This is the content of Theorem 76 in Howitt's Ph.D.\ thesis~\cite{H07}. The proof is based on martingale characterizations and a decomposition in terms of successive times of coalescence between paths.

\bigskip

\noindent
{\bf Acknowledgement} R.~Sun is supported by NUS grant R-146-000-253-114.

\bibliographystyle{abbrv}
\bibliography{bibweb}

\end{document}